\documentstyle[12pt]{article}

\newcommand{\const}{\mathop{\rm const}\limits}

\newcommand{\Law}{\mathop{\rm Law}\limits}

\begin{document}

\begin{center}

{\bf  CHARACTERIZATION OF QUASY-GAUSSIAN DISTRIBUTIONS. }\\

\vspace{4mm}

 $ {\bf E.Ostrovsky^a, \ \ L.Sirota^b,  \ \ A.Zeldin^c.} $ \\

\vspace{4mm}

$ ^a $ Corresponding Author. Department of Mathematics and computer science, \\
Bar-Ilan University, 84105, Ramat Gan, Israel.\\

\vspace{4mm}

E-mail:\ galo@list.ru \ eugostrovsky@list.ru\\

\vspace{4mm}

$ ^b $  Department of Mathematics and computer science. Bar-Ilan University, 84105,\\ Ramat Gan, Israel.\\

\vspace{4mm}

E-mail:\ sirota3@bezeqint.net\\

\vspace{4mm}

$ ^c $  Department of Mathematics and computer science. Bar-Ilan University, 84105, \\ Ramat Gan, Israel.\\

\vspace{4mm}

E-mail:\ anatolyz@moia.gov.il\\

\vspace{5mm}
                    {\bf Abstract.}\\

 \end{center}

 \vspace{4mm}

 A new characterization of the multivariate so-called "quasi-Gaussian distribution" (the authors dared to coin a new term) by means of independence their
Cartesian and polar coordinates proposed. \par
 The authors try to show that these distributions may essentially differ from the classical Gaussian distribution. Some properties
 of these distributions are studied: calculating of moments and of bilateral tail behavior.\par
Some potential applications of these distributions and mixtures of the distributions in demography, philology are discussed in the final section. \par

\vspace{4mm}


\vspace{4mm}


\vspace{4mm}

\section{ Statement of the problem. Notations. Definitions. }

\vspace{3mm}

{\bf A.}  There exist many characterizations of a two-dimensional, or, more generally, multidimensional Gaussian (normal)
distributions, with independent coordinates, for example, a characterization by means of independence of linear functionals  or through the
distribution of sums of coordinates, see the famous text - \cite{Feller1}, chapter III, p. 77 and chapter XV, p. 498-500;  by means of the
properties of conditional distributions \cite{Albajar1}, \cite{Kotlarski1}; a characterization by means of the properties of  order statistics
\cite{Jian1}; a characterization by means of some inequalities \cite{Bobkov1}, \cite{Kac1} etc., see also the reference therein.\par
 The classical monograph of Kagan, Linnik, Rao \cite{Kagan1} is completely devoted to the characterisation problems
in Mathematical Statistics.\par
 Usually, these characterizations are stable (robust)  \cite{Meshalkin1}, \cite{Zolotarev1}. \par

\vspace{3mm}

 Notice that the characterization by means of independence of linear functionals reduced as a rule to the following functional
 equation:

$$
f_1(x)  + f_2(y) = g_1(x+y) + g_2(x-y), \ x,y \in R, \ f_i, g_j: R \to R \eqno(1.1)
$$
 with unique solution up to multiplicative and additive constants  \cite{Feller1}, chapter XV, p. 498 - 500:

$$
f_i(z) = C_i \ z + D_i , \ g_j(z) = K_j z + L_j, \  C_i, D_i, K_j, L_j = \const,  \ i,j = 1,2; \eqno(1.2)
$$
evidently, $ C_1 = K_1 + K_2, \ C_2 = K_1 - K_2, \  D_1 + D_2 = L_1 + L_2.     $

\vspace{3mm}

{\it The immediate predecessors of the present paper are the articles  of Tamhankar \cite{Tamhankar1} and
Flusser \cite{Flusser1}, in which the two-dimensional Gaussian distribution with independent coordinates
is characterized by means of independence as well as in polar coordinates. } \par

\vspace{3mm}

  It was assumed in the aforementioned paper \cite{Tamhankar1}, \cite{Flusser1} that the considered r.v. $ \xi, \eta $
have continuous strictly positive densities $  f(x), \ g(y). $\par

\vspace{3mm}

{\bf  Our purpose is to show that if the technical conditions of positivity and continuity are not satisfied, and the
Cartesian (Descartes) and polar coordinates are pairwise independent, the random variable
may be not Gaussian; also we're trying to find the general form of these distributions. }\par

\vspace{3mm}

 We also offer applications for the obtained distributions. \par

\vspace{4mm}

{\bf B.} Let us consider the following examples.\\
{\bf Example 1.} Define the ordinary centered Gaussian non-degenerate density with variance $ \sigma^2, \ \sigma > 0 $  by

$$
f_{\sigma}(x) = (2 \pi)^{-1/2} \ \sigma^{-1} \ \exp \left(-x^2/(2 \sigma^2) \right), \ x \in R. \eqno(1.3)
$$
 Define for two numerical r.v.  $ (\xi, \eta) $  their  polar coordinates

$$
\rho = \rho(\xi, \eta) = \sqrt{\xi^2 + \eta^2}, \hspace{5mm} \theta = \theta(\xi, \eta) = \arctan(\eta/\xi). \eqno(1.4)
$$

 Suppose the r.v. $ \xi, \eta $ are centered independent identical distributed with normal density $ f_{\sigma}(x):
 \Law (\xi) = \Law(\eta) = N(0, \sigma^2).  $ Then the r.v. $ \rho(\xi, \eta) $ and $ \theta(\xi, \eta)  $  are also independent.
 Actually, we have for these values $ r = \const > 0, \phi \in (0, \pi/2) $ using polar coordinates

 $$
 {\bf P} (\rho <r, \ \theta < \phi) =  (2 \pi)^{-1} \ \sigma^{-2} \int_{ v < r } \int_{\beta < \phi}
 v \ \exp(-v^2/(2 \sigma^2) ) \ dv \ d \beta =
 $$

 $$
 (\phi/(2 \pi)) \times (1 - \exp(-r^2/(2 \sigma^2)). \eqno(1.5)
 $$
The case $ \theta > \pi/2 $ is considered analogously.\\
 Note in addition that if the random vector $ (\xi, \eta) $  has multivariate centered normal distribution with non-zero
correlation, then their polar coordinates are dependent. \par

\vspace{3mm}

{\it At first glance it may seem that this property characterizes unambiguously the two-dimensional normal distribution with
independent coordinates.}\par

\vspace{3mm}

 To ensure that the reverse, we consider the  second example (counterexample).\\

\vspace{3mm}

{\bf Example 2.} We denote as trivial for any measurable set $ A, \ A \subset R $ its indicator function by
$ I(A) = I_A(x):   $\par

$$
I_A(x) = 1, \ x \in A; \hspace{5mm} I_A(x) = 0, \ x \notin A.
$$

Introduce a {\it family} of functions

$$
\omega_{\alpha}(x) = \omega_{\alpha}(x; C_1, C_2) := C_1 \ |x|^{\alpha(1)} \ I_{(-\infty,0)}(x) + C_2 \ x^{\alpha(2)} \ I_{ (0,\infty)}(x),
$$

$$
x \in R, \ C_{1,2}= \const \ge 0, \ \alpha = \vec{\alpha}  = (\alpha(1), \alpha(2)),  \ \alpha(1), \alpha(2) = \const > -1, \eqno(1.6)
$$
so that $ \omega_{\alpha}(0) = 0,  $ and a family of a correspondent probability densities of a view

$$
g_{\alpha, \sigma}(x) = g_{\alpha, \sigma}(x; C_1, C_2)  \stackrel{def}{=} \omega_{\alpha}(x; C_1, C_2) \ f_{\sigma}(x). \eqno(1.7)
$$
 Since

$$
I_{\alpha(k)}(\sigma) := \int_0^{\infty} x^{\alpha(k)} \exp \left( -x^2/(2 \sigma^2)  \right) \ dx = 2^{(\alpha(k) - 1)/2  } \
\sigma^{(\alpha(k) + 1)} \ \Gamma((\alpha(k) + 1)/2),
$$
where $  \Gamma(\cdot) $ is ordinary Gamma function, there is the interrelation between the constants $ C_1, C_2  $ in (1.7):

$$
C_1 \ I_{\alpha(1)}(\sigma) + C_2 \ I_{\alpha(2)}(\sigma) = \sigma \ (2 \pi)^{1/2}, \eqno(1.8)
$$
has only one degree of freedom.  In particular, the constant $  C_1 $ may be equal  to zero;
in this case the r.v. $ \xi $ possess  only non-negative values.\par

\vspace{3mm}

{\it  We will denote by $  C_i, K_j $ some finite non-negative constants that are not  necessary to be
the same in different places. }\\

\vspace{3mm}

{\bf  C. \ Definition 1.1. } The distribution of a r.v. $  \xi   $ with density function of a view
$ x \to g_{\alpha, \sigma}(x - a; C_1, C_2), \ a = \const \in R $ is said  to be quasi-Gaussian or equally quasi-normal. Notation:

$$
\Law(\xi) = QN(a,\alpha,\sigma, C_1, C_2). \eqno(1.9)
$$

\vspace{3mm}

 The value $  "a" $ in (1.9) may be called  {\it  quasi-center  } by analogy with normal distribution  and value
 $ "\sigma" $ which is called {\it  quasi-standard } of the r.v. $  \xi. $ \\

\vspace{3mm}

{\bf D. Remark 1.1. Moments and tail behavior.} \\
 Denote $ \xi^p_+ = \xi^p \ I_{(0,\infty)}(\xi), \ \xi^p_- = \xi^p \ I_{(-\infty, 0)}(\xi). $  Evidently,
$$
 |\xi|^p =  \xi^p_+ + \xi^p_-, \ p \in R; \ \xi^n = |\xi|^n, \ n = 0, \pm 2, \pm 4, \ldots;
$$

$$
\xi^m = \xi^m_+ - \xi^m_-, \  m =  \pm 1, \pm 3, \pm 5 \ldots.
$$
 We have for the quasi-centered quasi-Gaussian distribution $ \Law(\xi) = QN(0, \vec{\alpha},\sigma, C_1, C_2): $

$$
{\bf E} \xi^p_+  = C_2 \ I_{\alpha(2) + p}(\sigma), \ \alpha(2) + p > -1; \ {\bf E} \xi^q_-  = C_1 \ I_{\alpha(1) + q}(\sigma),
\alpha(1) + q > -1.
$$

 Therefore

$$
{\bf E} |\xi| = C_2 \ I_{\alpha(2) + 1}(\sigma)+  C_1 \ I_{\alpha(1) + 1}(\sigma),
$$

$$
{\bf E} \xi = C_2 \ I_{\alpha(2) + 1}(\sigma) -  C_1 \ I_{\alpha(1) + 1}(\sigma),
$$

$$
{\bf E} \xi^2 = C_2 \ I_{\alpha(2) + 2}(\sigma)+  C_1 \ I_{\alpha(1) + 2}(\sigma).
$$

\vspace{3mm}
 Let us consider now the tail behavior of the r.v. $ \xi. $   We find using the Hospital's rule  as $ y \to +\infty: $

 $$
 {\bf P} (\xi  > y) \sim \frac{C_2 \ \sigma}{\sqrt{2 \pi}}  \ y^{\alpha(2) - 1} \ e^{- y^2/(2 \sigma^2)  },
 $$

 $$
 {\bf P} (\xi  < - y) \sim \frac{C_1 \ \sigma}{\sqrt{2 \pi}}  \ y^{\alpha(1) - 1} \ e^{- y^2/(2 \sigma^2)  }.
 $$

\vspace{4mm}

{\bf E. } \ It is easy to verify that if the r.v. $  (\xi, \eta) $ are independent and both have the quasi-Gaussian distribution
with parameters $  a = 0, \ b = 0 $ ("quasi-centered" case):

$$
\Law(\xi) = QN(0,\alpha,\sigma, C_1, C_2), \hspace{5mm} \Law(\eta) = QN(0,\beta,\sigma, C_3, C_4)\eqno(1.10)
$$
may be with different parameters $ \alpha \ne \beta, \ C_1 \ne C_3, C_2 \ne C_4  $ but with the same value of the standard
$  \sigma, \ \sigma > 0, $ then their polar coordinates  $ (\rho, \theta)  $ are also independent.\par

\vspace{3mm}

 We'll show further that the inverse proposition under simple natural conditions is true. \par

\vspace{3mm}

\section{Main result.}

\vspace{3mm}

{\bf Definition 2.1.}  The distribution on the real axis  of a random variable $  \xi  $  is said  to be regular
(at the origin), write:  $ \Law(\xi) \in Reg, $  if it has a density $ f(x) = f_{\xi}(x) $ such that there exists a constant $ \mu  $
for which

$$
0 < \inf_{x \in (-1,1) } \frac{f(x)}{|x|^{\mu}} \le  \sup_{x \in (-1,1) } \frac{f(x)}{|x|^{\mu}} < \infty. \eqno(2.0)
$$
 We write $ \mu = \deg{\xi} $   ("degree"). Of course, $ \mu > -1.  $ Evidently, the constants $ " -1", \ "1" $ in (2.0)
 may be replaced on other positive finite value. \\

\vspace{3mm}

{\bf Theorem.} {\it Let's suppose that the  random variables $  (\xi, \eta) $ are independent, both have the regular distribution,
 may be with different degree,  and their polar coordinates  $ (\rho, \theta)  $ are also independent. \par
Then  both the r.v. $  (\xi, \eta) $  have a  quasi-centered quasi-Gaussian distribution as in (1.10).  } \par

\vspace{3mm}

{\bf Remark 2.1.} We do not suppose here the continuity and strictly positivity of the densities of distributions of the r.v.
$  (\xi, \eta),  $ in contradiction to the articles \cite{Tamhankar1}, \cite{Flusser1}.\par

\vspace{3mm}

{\bf Proof.} \par
{\bf 1.}  Denote  by $ f(x), \ g(y) $ correspondingly the densities of the r.v. $ (\xi, \eta); $ then the common
density $ h(r, \phi) $ of their polar coordinates $ (\rho, \theta) $  may be written  as follows:

$$
h(r, \phi)  = \sqrt{x^2 + y^2} \ f(x) \ g(y). \eqno(2.1)
$$

 Evidently, the functions $ f,g,h $ are non - negative, measurable and integrable. \par
Because of the independence  the function $ h(r, \phi) $  is factorable:

$$
 h(r, \phi)  = h_1(r) \ h_2(\phi).
$$

We come to the following key functional equation of a view

$$
f(x) \ g(y) = M \left(\sqrt{x^2 + y^2} \right) \ L(y/x). \eqno(2.2)
$$
for some positive local integrable functions $ f,g,M,L.  $ \par

\vspace{3mm}

{\bf 2.} It suffices to consider the equation (2.2) only for positive values $ (x,y); $  other cases may be
investigated in a similar way. Then, we can rewrite  (2.2) after permutation $ (x,y) \to (y,x) $

$$
f(y) \ g(x) = M \left(\sqrt{x^2 + y^2} \right) \ L(x/y) \eqno(2.3)
$$
and we obtain after dividing (2.2) term by term by (2.3)

$$
\frac{f(x)}{g(x)} : \frac{f(y)}{g(y)} =    \frac{L(y/x)}{L(x/y)},   \eqno(2.4)
$$
or equally

$$
p(x) = p(y) \cdot \psi(x/y), \eqno(2.5)
$$
where we denote

$$
p(x) = \frac{f(x)}{g(x)}, \hspace{5mm} \psi(z) = \frac{L(1/z)}{L(z)}.
$$

\vspace{3mm}

{\bf 3.}  The next substitution $ p(x) = P( x) = \log p(x), \ \Psi(x) = \log \psi(\exp x) $
leads us to the equation

$$
P(x) = P(y) + \Psi(x-y), \eqno(2.6)
$$
with solution:  $  P(x) = C \cdot x. $ Therefore

$$
f(x) = x^{\gamma} \cdot  g(x), \ \gamma = \const, \ x > 0. \eqno(2.7)
$$
and we conclude substituting into (2.2)

$$
x^{\gamma} \cdot  g(x) \cdot g(y) = M \left(\sqrt{x^2 + y^2} \right) \ L(y/x). \eqno(2.8)
$$

 The example 2 shows that the constant $ \gamma $ may be non-zero. \par

\vspace{3mm}

{\bf 4.} The next substitution $  g(x) = x^{\lambda} w(x)  $ into (2.8) give us the following equation

$$
x^{\gamma + \lambda} \ y^{\lambda} \ w(x) \ w(y) = M \left(\sqrt{x^2 + y^2} \right) \ L(y/x). \eqno(2.9)
$$
 If in (2.9) we choose  $ \gamma + \lambda = - \lambda, $ i.e. $ \lambda = -\gamma/2,  $ we get
instead  (2.8)  the equation

$$
  w(x) \cdot w(y) = M \left(\sqrt{x^2 + y^2} \right) \ L_1(y/x) \eqno(2.10)
$$
with another in general case homogeneous function $ L_1(y/x). $ \par
 Thus, we reduce the general case to the possibility when the r.v. $  \xi, \eta $ are in addition identically distributed.\par

\vspace{3mm}

{\bf 5.} The another change of variables  $ w(x) = x^{\nu} \ w(x), \nu = \const $ leads us to the functional equation  of a view

$$
  v(x) \cdot v(y) = M_1 \left(\sqrt{x^2 + y^2} \right) \ L_2(y/x), \eqno(2.11)
$$
however due to the proper selection of parameter $  \nu  $ it can be achieved that

$$
0 < \inf_{x \in (0,1) } v(x) \le  \sup_{x \in (0,1) } v(x) < \infty. \eqno(2.12)
$$
 The solution of (2.11) having a probabilistic sense under the condition (2.12) looks like
 \cite{Tamhankar1}, \cite{Flusser1}:

 $$
v(x) = C_6 \cdot \exp(- C_7 x^2), \ C_6, C_7 = \const > 0. \eqno(2.13)
 $$
We should return to the source functions $  f(x), \ g(y). $ \par

\vspace{3mm}

{\bf 6.}  Suppose in addition that all the functions in (2.11) are positive and continuous differentiable almost everywhere. Then

$$
  v_1(x) + v_1(y) = M_2 \left(\sqrt{x^2 + y^2} \right) + L_3(y/x), \eqno(2.14)
$$
where

$$
v_1(x) = \log v(x), \ M_2(x) = \log M_1(x), \ L_3(x) = \log L_2(x).
$$

 Denote $ z(x) = v_1(x^2), $ then

$$
z(x) + z(y) = M_2(x + y) + L_4(y/x).  \eqno(2.15)
$$

Without loss of generality we can adopt $  z(0) = M_2(0) = L_4(0) = 0, \ 0/0 = 0, \ L_4(0/0) = L_4(0) = 0.  $  Putting in (2.15)
 the value $  y = 0,  $ we deduce $  M_2(x) = z(x); $  therefore

 $$
 z(x+y) - z(x) = z(y) + L_5(y/x). \eqno(2.16)
 $$
  We choose in (2.16) $ y = \Delta x, \ \Delta x \to 0:  $

$$
z(x + \Delta x) - z(x) = C_1 \ \Delta x  + L_5 (\Delta x/x) = C_2 \Delta x + o(\Delta x), \ dz /dx = C_2 = \const.
$$
  Thus, $ z(x) = C_2 x + C_3. $ The proof of the theorem is completed. \par

  \vspace{4mm}

\section{Concluding remarks.}

\vspace{3mm}

{\bf A.  Multidimensional case.}\par

 It is possible to generalize our theorem on the multidimensional case. Actually, let us consider the random vector
$ \xi = \vec{\xi} =  (\xi_1, \xi_2, \ldots, \xi_d) $  with the density

$$
f_{\xi}(x_1, x_2, \ldots, x_d) = G( x_1, x_2, \ldots, x_d; \vec{\alpha}, \vec{\sigma}, \vec{ C_1 },  \vec{C_2  } ) \stackrel{def}{=}
$$

$$
\prod_{j=1}^d  g_{\alpha_j, \sigma_l}(x_j; C_1^{(j)}, C_2^{(j)}), \eqno(3.1)
$$
where $ \alpha_j > -1, \ \sigma_j = \const  > 0, \ C_i^{(j)}  = \const \ge 0, $

$$
C_1^{(j)} \ I_{\alpha^{(j)}(1)}(\sigma_j)  + C_2^{(j)} \ I_{\alpha^{(j)}(2)}(\sigma_j) = \sigma_j \ (2 \pi)^{1/2}, \eqno(3.2)
$$

\vspace{3mm}

{\bf Proposition 3.1.} \par

{\it Assume that all the standards $ \sigma_j = \sigma $ do not dependent on the number $  j. $ Then
the (Cartesian) coordinates of the vector $  \vec{\xi}, $ i.e. the random variables $  \{  \xi_j \} $ are common independent
and so are their polar coordinates. \par
 The contrary is also true: if
the Cartesian and polar coordinates of the vector $  \vec{\xi} $ are common independent and the random variables $  \{  \xi_j \} $
are regular distributed, then it density has a view (3.1), with the same standards } $ \sigma. $ \par

\vspace{3mm}

{\bf B. Weight quasi-Gaussian distributions. } \par

 \vspace{3mm}

 Let $ W_k, \ k = 1,2, \ldots, N $ be positive numbers (weights) such that $ \sum_{k=1}^N W_k = 1. $ We define the {\it weight}
or {\it mixed} quasi-Gaussian distribution by means of multivariate density of a view

$$
G^{(W)}( x_1, x_2, \ldots, x_d) = G^{(W)} \left( x_1, x_2, \ldots, x_d; \{ a_j^{(k)} \}, \{ \alpha_j^{(k)} \}, \{ \sigma_j^{(k)} \}, \{ C_1^{(k)} \}\right) \stackrel{def}{=}
$$

$$
\sum_{k=1}^N  W_k \ G \left( x_1-a_1^{(k)}, x_2-a_2^{(k)}, \ldots, x_d-a_d^{(k)}; \vec{\alpha}^{(k)}, \ \{ \sigma_j^{(k)} \}, \ \vec{ C_1 }^{(k)} \right).
\eqno(3.3)
$$

\vspace{3mm}

 More general view of similar distribution has a discrete component:

 $$
 G_0^{(W)}( \vec{x}) := W_0 \delta(\vec{x} - \vec{a_0}) +
 $$

 $$
 \sum_{k=1}^N  W_k \ G \left( x_1-a_1^{(k)}, x_2-a_2^{(k)}, \ldots, x_d-a_d^{(k)}; \vec{\alpha}^{(k)}, \vec{\sigma}^{(k)}, \ \vec{ C_1 }^{(k)} \right), \eqno(3.4)
 $$

$$
W_0, W_1, \ldots, W_N > 0, \ \sum_{k=0}^N W_K = 1,
$$

$ \delta(\vec{x}) $ is the classical Dirac delta function; so that

$$
{\bf P} (\vec{\xi} =  \vec{a_0}) = W_0 > 0.
$$

{\bf C. Possible applications. } \par

\vspace{3mm}

 The classic application  of  ordinary characterization theorem in  Statistical Physics is described, e.g. in the book \cite{Kagan1},
 chapter 3, section 3.5. \par

\vspace{3mm}

 By our opinion,  the weighted quasi-Gaussian distribution may be used as a first approximation in the {\it demography,} where  $ (\xi_1, \xi_2) $
may be  coordinates of a place of residence of a  random  person, and $ W_k $ is a share of $ k^{th} $ city in the general population of some
country. \par

\vspace{3mm}

 Another possible application is in the realm of {\it philology,}  where $ \vec{\xi} = (\xi_1, \xi_2, \ldots, \xi_d) $  may denote the parameters of a certain word, for instance, phonetic and/or semantic values,
  and $ N; \{  a_j^{(k)} \} $ may be estimated by means of cluster  analysis and by methods of
parametrical statistics: maximum likelihood, chi-square test etc.\par
Thus a lexical unit can be analyzed either as an independent entity by itself, or being compared to (a). other units from the same language/other related or unrelated languages of the same epoch (synchronic view); (b) historically related words of the similar semantic layer (diachronic view). The described method might prove to be very useful also in lexicostatistics.
 The quasi - centers $ \{  a_j^{(k)} \} $ may be interpreted as a coordinates of fundamental human notions:  food, policy, medicine, economic etc. \par
 The classification based on the  mixed quasi-Gaussian distribution
may be useful in learning a foreign language.\par

\vspace{3mm}

Needless to say, this approach requires an experimental verification.\par

\end{document}